\documentclass[12pt,a4paper]{article}
\usepackage{amssymb,amsthm,amsmath}

\allowdisplaybreaks[4]
\newenvironment{prova}{\begin{proof}[Proof:]\parindent=0in~\\ }{\end{proof}}

\newtheorem{thm}{Theorem}[section]

\newtheorem{prop}[thm]{Proposition}

\newtheorem{df}[thm]{Definition}

\newcommand{\bn}[2]{\Big(\textrm{\footnotesize$\!\!\begin{array}{c}#1 \\ #2\end{array}\!\!$}\Big)}
\newcommand{\A}{\mathcal{A}}
\newcommand{\B}{\mathcal{B}}
\newcommand{\K}{\mathcal{K}}
\newcommand{\D}{\mathcal{D}}
\newcommand{\HH}{\mathcal{H}}
\newcommand{\HHH}{\hat{\mathcal{H}}}
\newcommand{\T}{\mathcal{T}}
\newcommand{\rd}{\mathcal{S}} % Rapid decay sequences
\newcommand{\C}{\mathbb{C}}
\newcommand{\R}{\mathbb{R}}
\newcommand{\N}{\mathbb{Z}_+}
\newcommand{\NN}{\mathbb{N}}
\newcommand{\Z}{\mathbb{Z}}
\newcommand{\ket}[1]{\left|#1\right>}
\newcommand{\inner}[1]{\left<#1\right>}
\newcommand{\tr}{\mathrm{Trace}}
\newcommand{\de}{\mathrm{d}}
\newcommand{\ma}[1]{\left(\!\begin{array}{cc}#1\end{array}\!\right)}

\begin{document}
\title{~\\ Local Index Formula on the Equatorial Podle\'s Sphere}
\date{~\\}
\author{~\\
Francesco D'Andrea and Ludwik D\c{a}browski
\\[30pt] 
{\small Scuola Internazionale Superiore di Studi Avanzati,}\\[5pt]
{\small Via Beirut 2-4, I-34014, Trieste, Italy} \\[8pt]
}
\maketitle
                                                                               
\begin{abstract}
\noindent 
We discuss spectral properties of the equatorial Podle\'s sphere $S_q^2$.
As a preparation we also study the `degenerate' (i.e.~$q=0$) case
(related to the quantum disk).
Over $S_q^2$ we consider two different spectral triples: 
one related to the Fock representation of the Toeplitz algebra
and the isopectral one given in~\cite{DLPS}. 
After the identification of the smooth pre-$C^*$-algebra we compute
the dimension spectrum and residues. We check the nontriviality 
of the (noncommutative) Chern character of the associated Fredholm modules
by computing the pairing with the fundamental projector of the $C^*$-algebra
(the nontrivial generator of the $K_0$-group) as well as the pairing
with the $q$-analogue of the Bott projector. Finally, we show that the
local index formula is trivially satisfied.
\end{abstract}

\vfill\eject

%%% ======================================================================
\section{Introduction.}
The noncommutative differential (or \emph{spectral}) geometry 
of quantum groups and $q$-deformed spaces has been recently
intensively studied (see e.g.~\cite{Dab} for review).
In particular the explicit computation of the local index formula for the
total characteristic class, i.e.~the Connes-Chern character, 
has been worked out on the quantum group $SU_q(2)$ 
in~\cite{Con1} for the `singular' spectral triple of~\cite{CP} 
and in~\cite{vSDLSV} for the %even
spectral triple of~\cite{DLSvSV}, with most of the results coinciding.

In this paper we present a systematic discussion of analogous 
spectral properties of another quantum space: 
the equatorial Podle\'s sphere $S_q^2$, originally defined for $q\neq 0$
as a homogeneous $SU_q(2)$-space.
We also study separately the `degenerate' (i.e.~$q=0$) case,
that has a perfect meaning (though not as a homogeneous space).

Our main task is the analysis of the local index formula for the spectral 
triple on $S_q^2$ constructed in~\cite{DLPS}, which anticipated 
some of the interesting properties of that in~\cite{DLSvSV}.
We analyze also another spectral triple, related to the Fock representation 
of the Toeplitz algebra.

For that purpose it is convenient to study first the `degenerate' 
(i.e.~$q=0$) case (related to the quantum disk),
on which we consider two different spectral triples.
Using these results, on $S_q^2$ we analyze two spectral triples:
the one related to the Fock representation of the Toeplitz algebra
(in Section~\ref{sec:mu})
and the isopectral one given in~\cite{DLPS} (in Section~\ref{sec:pi}). 
After the identification of the pre-$C^*$-algebra of `smooth' elements, 
we compute the dimension spectrum and residues.
We check the nontriviality of the (noncommutative) Chern character
of the associated Fredholm modules by computing the pairing
with the fundamental projector of the $C^*$-algebra
(the nontrivial generator of the $K_0$-group) as well as the pairing
with the $q$-analogue of the Bott projector. Finally, we explicitely 
verify that the local index formula is trivially satisfied.

The relevance of such explicit calculations and results 
stems from their relative scarceness in the literature
for interesting noncommutative spaces.

In the following we use the notation $\Z_+=\{1,2,3,\ldots\}$ and $\NN=\Z_+\cup\{0\}$.

%%% ======================================================================

\section{Preliminaries about the equatorial Podle\'s sphere.}\label{sec1}

We use, with minor changes, notation of~\cite{DLPS}. For $0<q<1$ the $*$-algebra
(of polynomials) $\A(S_q^2)$ on the equatorial Podle\'s sphere is generated by 
$a$, $a^*$ and $b=b^*$ with relations\footnote{The original presentation of
Podle\'s~\cite[eq.~(7b)]{Pod} corresponds to $A=q^{-2}b$, $B=a^*$ and $\mu=q$
(for $q<1$).}
\begin{equation}\label{eq:def}
ba=q^2ab\ , \qquad\quad a^*a+b^2=1 \ , \qquad\quad q^4 aa^*+b^2=q^4\ .
\end{equation}
Its linear basis can be taken as $\{a^nb^m,\;(a^*)^{n+1}b^m,\;n,m\in\,\NN\}$.

We denote by $C(S^2_q)$ the universal $C^*$-algebra of $\A(S^2_q)$.
$S^2_q$ is known to be an embeddable $SU_q(2)$-homogeneous space, and carries
a strongly continuous action of $S^1$. Being embeddable, the coaction
of $SU_q(2)$ it carries is the restriction of the coproduct of $SU_q(2)$.
One easily verifies that the subalgebra $\A(S^2_q)\subset\A(SU_q(2))$ is invariant
with respect to the standard left action~of~$U_q(su(2))$.

The relations (\ref{eq:def}) make sense for any $q\in\R$ but 
the map 
$a\mapsto a^*$, $b\mapsto q^{-2}b$ extends to an isomorphism between 
$\A(S^2_q)$ and $\A(S^2_{q^{-1}})$, while one has trivially 
$\A(S^2_q)\simeq \A(S^2_{-q})$. 
So, without loss of generality, one can restrict to $0\leq q\leq 1$.
The two limiting cases are as follows:
\begin{enumerate}
\item If $q=1$, this is just the `polynomial' algebra on the (commutative) sphere $S^2$,
while its $C^*$-algebra closure corresponds to continuous functions on $S^2$.
\item If $q=0$, we have $b=0$ and denoting $a=w$ we see that
$\A(S^2_0)$ is the algebra of polynomial functions in the elements $w$ and
$w^*$, with relation $\,w^*w=1$\ .
\end{enumerate}
The associated universal $C^*$-algebra $C(S^2_0)$, being 
generated by one partial isometry, is known to be isometrically
$*$-isomorphic to the algebra of Toeplitz operators $\T$~\cite{Cob}.
We can also think of $C(S^2_0)\simeq\T$ as the algebra $C(\D_q)$ 
of continuous `functions' on a non-commutative disk $\D_q$, 
due to the short exact sequence\footnote{Here $\sigma$ is the
symbol map, sending the unilateral shift $w$ to $e^{i\theta}$ and extended to a
$C^*$-algebra morphism.}
\begin{equation}\label{eq:uno}
0\to\K\to\T\stackrel{\sigma}{\to} C(S^1)\to 0\ ,
\end{equation}
interpreted as the non-commutative analogue of the sequence
\begin{equation*}
0\to C_0(\textrm{open disk})\to C(\textrm{closed disk})\to C(S^1)\to 0\ .
\end{equation*}
In eq.~(\ref{eq:uno}), we think of $S^1$ as the boundary of the non-commutative
disk, $\sigma$ as the evaluation on the boundary and the compact operators $\K$
as continuous `functions' on $\D_q$ vanishing on the boundary.

In the sequel, unless stated otherwise we shall require that $0<q<1$.\\
From the last defining relation in (\ref{eq:def}) it easily follows that 
in all representations the operator norm of $q^{-2}b$ is $\leq 1$. 
Hence the $C^*$-norm of $b$ satisfies $||b||<1$.
Therefore, $1-b^2=a^*a$ is invertible in the $C^*$-algebra and the element
\begin{equation}\label{eq:p}
p=1-a(a^*a)^{-1}a^*\in C(S^2_q)
\end{equation}
is a projector (i.e.~$p=p^*=p^2$). In Appendix~\ref{app:B} we show that the projective
left $C(S^2_q)$-module $L^2(S^2_q)p$ is equivalent to the (graded) faithful representation
given in~\cite{Pod} $\mu:=\mu_+\oplus\mu_-:C(S^2_q)\to\B(\ell^2(\N)\oplus\ell^2(\N))$, where
\begin{equation}\label{eq:qdisk}
\mu_{\pm}(a)\ket{n}=\sqrt{1-q^{4n}}\ket{n+1}\;,\qquad\quad
\mu_{\pm}(b)\ket{n}=\pm q^{2n}\ket{n}\;,
\end{equation}
with $\ket{n}$ being the canonical orthonormal basis of $\ell^2(\N)$.

The representations $\mu_{\pm}$ are irreducible but not faithful;
$ \mu_+(b)$ has positive spectrum, while $ \mu_-(b)$ has negative spectrum.
For this reason, we interpret $C(S^2_q)/\ker\mu_{\pm}$ as a couple of closed
(noncommutative) hemispheres composing the sphere.

We define $\A(\D_q):= \mu_{+}(\A(S^2_q))=\mu_{-}(\A(S^2_q))$, i.e.~as the polynomial
$*$-algebra generated by $\mu_{\pm}(a)$ and $\mu_{\pm}(b)$, given by (\ref{eq:qdisk}).
(These are the noncommutative disks considered in~\cite{Con1, vSDLSV}).\\
Notice that at the $C^*$-algebra level
\begin{equation*}
C(S^2_q)/\ker\mu_{\pm}\simeq C(\D_q)=\T\simeq C(S^2_0)\;,
\end{equation*}
even though $\A(\D_q)\neq\A(S^2_0)$.

\pagebreak

Identifying $C(S^2_q)$ with $\mu (C(S^2_q))$, in~\cite{Sheu} it was shown that
the map $f\mapsto \left( \mu_+(f), \mu_-(f)\right) $ gives the isomorphism 
\begin{equation}\label{eq:qD}
C(S^2_q)\simeq\{(x,y)\in C(\D_q)\oplus C(\D_q)\,|\,\sigma(x)=\sigma(y)\}\;.
\end{equation}
Thus, a `function' on the equatorial Podle\'s sphere is given by a couple 
of `functions' on the two hemispheres which coincide on the boundary \
($\sigma(x)=\sigma(y)$). This allows an interpretation of $S^2_q$ 
as a couple of noncommutative disks glued along the common boundary $S^1$.\\
In the degenerate case $q=0$, one can think that the two disks
`collapse' one over the other, as~$C(S^2_0)\simeq C(\D_q)$.

Notice that $C(S^2_q)\simeq C(S^2_{q'})$ for all $0<q,q'<1$. Furthermore
the representation $\mu$ allows us to construct a regular even spectral triple
over $\ell^2(\N)\oplus\ell^2(\N)$ with Dirac operator $(N\oplus N)F$, 
where $F$ flips the two subspaces and 
\begin{equation*}
N\ket{n}=n\ket{n}
\end{equation*}
is the `number' operator. 
This triple (modulo an irrelevant shift of $N$ by $1$) was studied in~\cite{CP1}. 

The kernel of the projection $C(S^2_q)\to\T$, $(x,y)\mapsto x$, is $\K$.
The map $\T\to C(S^2_q)$, $x\mapsto (x,x)$, proves that the following
sequence is split exact:
\begin{equation*}
0\to\K\to C(S^2_q)\leftrightarrow\T\to 0
\end{equation*}
and then $K_0(S^2_q)=\Z\oplus\Z$ and $K_1(S^2_q)=0$~\cite[cor.~8.2.2]{Wegge}.

The $C^*$-algebra morphism $\rho :C(S^2_q)\to C(S^1)$, $\rho(x,y):=\sigma(x)=\sigma(y)$,
gives rise to the exact sequence
\begin{equation*}
0\to\K\oplus\K\to C(S^2_q)\stackrel{\rho}{\to} C(S^1)\to 0\;.
\end{equation*}
The map $\rho$ is just the $C^*$-algebra morphism that extends the
map $a\to e^{i\theta}$, $b\to 0$.

At the level of polynomial algebras, we have the isomorphism:
\begin{equation*}
\A(S^2_q)\simeq\{(x,y)\in \A(\D_q)\oplus\A(\D_q)\,|\,\sigma(x)=\sigma(y)\}\;.
\end{equation*}

In Section.~\ref{sec:A} a certain (Fr\'echet) pre-$C^*$-algebra $A^\infty\subset\T$
is introduced. It fits the exact sequence
\begin{equation*}
0\to\rd\to A^\infty\to C^\infty(S^1)\to 0\;,
\end{equation*}
where $\rd$ is isomorphic to the pre-$C^*$-algebra of rapid decay matrices on $\ell^2(\NN)$.
There it is shown that $A^\infty$ (the definition of which does not depend on $q$) 
contains both $\A(S^2_0)$ and $\A(\D_q)$ for any $q$, so it can be interpreted as
the algebra of smooth elements (and used to construct a regular spectral triple)
on both $S^2_0$ and $\D_q$.

Since $S^2_q$ is obtained by gluing two copies of $\D_q$ along the boundary,
we define the smooth `functions' over $S^2_q$ as
\begin{equation}\label{eq:Cinf}
C^\infty(S^2_q):=\{(x,y)\in A^\infty\oplus A^\infty\,|\,\sigma(x)=
\sigma(y)\}\subset C(S^2_q)\;.
\end{equation}
It is a pre-$C^*$-algebra independent on $q$, and fits the exact sequence
\begin{equation*}
0\to\rd\oplus\rd\to C^\infty(S^2_q)\to C^\infty(S^1)\to 0\;.
\end{equation*}

Let us describe now the spin representation of $C(S^2_q)$, as defined
in~\cite{DLPS}. For that we shall use the Hilbert space isomorphic to 
classical $L^2$-spinors over the round sphere $S^2$, 
 $L^2(S^2)\otimes\C^2 =\HH_+\oplus \HH_-$, 
with orthonormal basis of `spinor harmonics' $\ket{l,m}_{\pm}$ 
labeled by $l\in\NN+\tfrac{1}{2}$ and $m=-l,-l+1,\ldots,l$. 
In this basis, the Dirac operator (over $S^2$)
is $D\ket{l,m}_{\pm}=(l+\tfrac{1}{2})\ket{l,m}_{\mp}$.\\
This Hilbert space can be also viewed as a module over $\A(S^2_q)$ for any $0<q<1$, 
that naturally extends (by continuity) to a bounded representation of the $C^*$-algebra.
The \emph{chiral} representations of~\cite{DLPS} are the faithful and irreducible
representations $\pi_{\pm}:C(S^2_q)\to\B(\HH_{\pm})$ defined by
\label{eq:Drep}
\begin{align*}
\pi_{\pm}(a)\ket{l,m}_{\pm} := &\;q^{m-l-\tfrac{1}{2}}\frac{
        \sqrt{[l+m+1][l+m+2]}}{[2l+2]}\ket{l+1,m+1}_{\pm} \\
      & -q^{m+l+\tfrac{1}{2}}\frac{\sqrt{[l-m-1][l-m]}}{[2l]}\ket{l-1,m+1}_{\pm} \\
      & \pm\frac{(1+q^2)q^{m-\tfrac{1}{2}}}{[2l][2l+2]}\sqrt{[l+m+1][l-m]}
        \ket{l,m+1}_{\pm}\;,
\\ %
\pi_{\pm}(b)\ket{l,m}_{\pm} := &\,-q^{m+1}\frac{\sqrt{[l+m+1][l-m+1]}}{[2l+2]}\ket{l+1,m}_{\pm} \\
      & -q^{m+1}\frac{\sqrt{[l+m][l-m]}}{[2l]}\ket{l-1,m}_{\pm} \\
      & \pm \frac{[l-m+1][l+m]-q^2[l-m][l+m+1]}{[2l][2l+2]}\ket{l,m}_{\pm}\;,
\end{align*}
where $[x]=(q^x-q^{-x})/(q-q^{-1})$ is the $q$-analogue of $x\in\C$.
The \emph{spin representation} $\pi$ over $\HH:=\HH_+\oplus\HH_-$ is the direct sum
$\pi:=\pi_+\oplus\pi_-$.\\
When $q=1$, $\pi$ is just the representation of $C(S^2_q)$ defined by left
multiplication on sections of the spin bundle over $S^2$.

We verify in Section~\ref{sec:mu} that $(C^\infty(S^2_q),\HH,D)$ is an (even) regular
spectral triple (isospectral, since $\HH$ and $D$ are the classical ones).
Isospectrality also means that the noncommutative Sobolev spaces (defined as domains
of $|D|^s$) are the ordinary ones over $S^2$; the same for the smooth domain of $D$
and for smoothing operators.

Even though originally defined for $q\neq 0$, 
the representations $\mu_\pm$ and $\pi_\pm$ make sense also for $q=0$.
In this case, both $\mu_{\pm}$ reduces to the (irreducible faithful) 
Fock representation of the Toeplitz algebra, with $w$ 
acting as the unilateral shift on $\ell^2(\N)$ 
\begin{equation}
\label{mu0}
\mu_{\pm}(w)\ket{n}=\ket{n+1}\;.
\end{equation}
Instead the representations $\pi_\pm$ of~\cite{DLPS} become
\begin{equation}
\label{pi0}
\pi_{\pm}(w)\ket{l,m}=\ket{l+1,m+1}\;,
\end{equation}
which when projected to a subspace with $l-m=k$ fixed, 
are equivalent to the Fock one.\\
Since in this degenerate $q=0$ case, $\mu_+=\mu_-$ and $\pi_+=\pi_-$,
the sign of the Dirac operator commutes with the algebra and is irrelevant.
For this reason in the next section we %shall
consider only one copy of (chiral) representation and the absolute value
of the Dirac operator.

%%% ======================================================================
\section{Spectral triples over the quantum disk}\label{sec:q0}

We shall use the notation of Section~\ref{sec1}.

%%% ======================================================================
\subsection{Description of the algebra of ``smooth'' elements}\label{sec:A}
We define $A^\infty\subset\T$ as the linear span of the elements
\begin{equation*}
f=\sum_{n\in\NN}(f_nw^n+f_{-n-1}(w^*)^{n+1})+\sum_{j,k\in\NN}f_{jk}w^j
(1-ww^*)(w^*)^k\;,\qquad\{f_n\}\in\rd(\Z)\,,\;\{f_{jk}\}\in\rd\,,
\end{equation*}

\vspace{-5pt}

\noindent%
where $\rd(\Z)$ indicates rapid decay sequences and $\rd$
rapid decay matrices on $\ell^2(\NN)$.

By direct calculation one checks that $A^\infty$ is a $*$-algebra (associative,
with unit) and the map $\sigma:A^\infty\to C^\infty(S^1)$, $f\mapsto
\sigma(f)$, $\sigma(f)(\theta):=\sum_{n\in\Z}f_ne^{in\theta}$, is a surjective $*$-algebra
morphism (this follows from the simple observation that, via Fourier series,
$\rd(\Z)$ equipped with convolution product is isomorphic to $C^\infty(S^1)$).
We used the symbol $\sigma$ because it is just the restriction to $A^\infty$
of the symbol map in eq.~(\ref{eq:uno}), i.e.~the $C^*$-algebra morphism
defined by $\sigma(w)(\theta):=e^{i\theta}$.

Furthermore, $\ker\sigma$ is a two-sided $*$-ideal in $A^\infty$ isomorphic
to $\rd$. At the level of abstract algebras, this follows from the equality
\begin{equation*}
w^j(1-ww^*)(w^*)^kw^{j'}(1-ww^*)(w^*)^{k'}=\delta_{j'k}w^j(1-ww^*)(w^*)^{k'}
\end{equation*}
and becomes more evident in the Fock representation, where
\begin{equation*}
f\ket{k}=\sum\nolimits_{j\in\N}f_{j-1,k-1}\ket{j}\;,
\qquad\quad
\forall\;f=\sum\nolimits_{j,k\in\NN}f_{jk}w^j(1-ww^*)(w^*)^k\in\ker\sigma\;.
\end{equation*}
Then, we have the short exact sequence anticipated in Section \ref{sec1}:
\begin{equation}\label{eq:seq}
0\to\rd\to A^\infty\stackrel{\sigma}{\to} C^\infty(S^1)\to 0\;.
\end{equation}
Notice that the action of $S^1$ on $\A(S^2_0)$, given by
$w\mapsto e^{i\theta}w$, defines a one-parameter group of automorphism
implemented on $\ell^2(\N)$ by the unitary operators $e^{i\theta N}$,
i.e.~for each $x\in \A(S^2_0)$ the action is $x\mapsto e^{i\theta N}xe^{-i\theta N}$.
Being implemented by unitary operators, it extends to a (strongly continuous) action
of $S^1$ on the $C^*$-algebra $\T$, and then on $A^\infty$.
The map $\sigma$ commutes with this action (it is an $S^1$-module morphism),
and elements of $A^\infty$ are smooth for the action of~$S^1$.

We conclude with the statement,
\begin{prop}
$A^\infty$ is a Fr\'echet pre-$C^*$-algebra.
\end{prop}
This follows form the fact that as a vector space $\A^\infty$ 
is the cartesian product of two Fr\'echet spaces, hence it is a Fr\'echet space too. 
To prove that it is a pre-$C^*$-algebra
it is sufficient to show that it contains the inverse of each element $f\in A^\infty$,
whenever $f$ is invertible in $\T$. 
To reach this conclusion, one can easily adapt the proof
of proposition 1 in~\cite{Con1} to the present case.

Clearly, $A^\infty$ contains $\A(S^2_0)$ and is dense in $\T$. 
But $A^\infty$ contains also the algebra $\A(\D_q)$. Indeed,
\begin{equation*}
q^{2N}=\sum\nolimits_{j,k\in\NN}\delta_{jk}q^{2k}w^{k+1}(1-ww^*)(w^*)^{k+1}
\end{equation*}
is an element of $A^\infty$. 
Now, $A^\infty$ being closed under holomorphic functional calculus, 
also $\sqrt{1-q^{4N}}\in A^\infty$.
Notice that $\A(\D_q)$ is generated by $q^{2N}$ and $w\sqrt{1-q^{4N}}$
(c.f.~Equation (\ref{eq:qdisk})).
So, $A^\infty$ contains the generators of $\A(\D_q)$ and then 
all the polynomial algebra.\\
For this reason, we can identify $A^\infty$ with smooth ``functions'' over the
quantum disk.
(Notice that the isomorphism $C(S^2_0)=\T\stackrel{\sim}{\to}C(\D_q)$ 
is ``smooth'' but not polynomial).

%%% ======================================================================
\subsection{Description of spectral triples for $A^\infty$}

We consider two spectral triples for $A^\infty$.
The first one is associated with the natural
representation (\ref{mu0}) of $A^\infty$ on $\ell^2(\N)$,
with $N$ taken as Dirac-type operator.
The second one is associated with the representation (\ref{pi0})
(i.e.~the $q\to 0$ limit of the isospectral representation of~\cite{DLPS}) 
together with $|D|$ (the absolute value of the classical Dirac operator $D$).

If $\ket{n}$ is the canonical basis of $\ell^2(\N)$ and
$\ket{l,m}$ the orthonormal basis of $\HH_+\simeq\bigoplus_{l+1/2\in\N}
\!\!\C^{2l+1}$ considered in~\cite{DLPS}, we are thus dealing respectively 
with the following representations
\begin{equation*}
w\ket{n}=\ket{n+1}\:,\qquad w\ket{l,m}=\ket{l+1,m+1}\;,\qquad\;n,l+\tfrac{1}{2}
\in\N,\;m=-l,-l+1,\ldots,l
\end{equation*}
and with the following `Dirac' operators
\begin{equation*}
N\ket{n}=n\ket{n}\:,\qquad\quad |D|\ket{l,m}=(l+\tfrac{1}{2})\ket{l,m}
\end{equation*}
(we identify $A^\infty$ with its representation).
Since both these operators are positive, the associated index map is trivial.

From the equations $[N,w]=w$ and $[|D|,w]=w$ we deduce that $A^\infty$
is invariant with respect to ~both the derivations $[N,\,.\,]$ and $[|D|,\,.\,]$, 
hence it is in their smooth domain. We have thus proved,
\begin{prop}
$(A^\infty,\ell^2(\N),N)$ and $(A^\infty,\HH_+,|D|)$ are regular
spectral triples.
\end{prop}
Moreover
\begin{prop}\label{prop:ker}
$\tr_{\ell^2(\N)}(fN^{-s})$ and $\tr_{\HH_+}(f|D|^{-s})$ are holomorphic functions
on $\C$, for all $f\in \rd = \ker\sigma$.
\end{prop}
\begin{prova}
If $f=\sum_{j,k\in\NN}f_{jk}w^j(1-ww^*)(w^*)^k$ is a generic element 
of $\ker\sigma$, then $f_{jj}\in\rd(\NN)$.

In both the triples, $w^j(1-ww^*)(w^*)^k$ is off-diagonal if $j\neq k$, 
while if $j=k$ it is a projector on $\ket{k+1}$, respectively on
$\C$-span of $\{\ket{k+\tfrac{1}{2},m}|~m=-k-\tfrac{1}{2}, \dots k+\tfrac{1}{2}\} $.
Thus we have the equalities
\begin{align*}
\tr_{\ell^2(\N)}(fN^{-s}) &=\sum\nolimits_{n\in\N}n^{-s}f_{n-1,n-1}\;, \\
\tr_{\HH_+}(f|D|^{-s}) &=\sum\nolimits_{l+\tfrac{1}{2}\in\N}(2l+1)(l+\tfrac{1}{2})^{-s}
f_{l-\tfrac{1}{2},l-\tfrac{1}{2}}
\end{align*}
and the series converge to holomorphic functions on $\C$ (Weierstrass theorem).
\end{prova}

\begin{prop}\label{prop:proof}
The dimension spectrum is $\{1\}$ for the first triple and $\{2\}$ for the second.
The residues are
\begin{equation*}
\mathrm{Res}_{s=1}\tr_{\ell^2(\N)}(fN^{-s})=\frac{1}{2\pi}\int_{S^1}\sigma(f)\de\theta
\;,\qquad
\mathrm{Res}_{s=2}\tr_{\HH_+}(f|D|^{-s})=\frac{1}{\pi}\int_{S^1}\sigma(f)\de\theta
\end{equation*}
for all $f\in A^\infty$, where $\sigma$ is the map in (\ref{eq:seq}).
\end{prop}
\begin{prova}
The algebra generated by $A^\infty$ and the commutators with the Dirac operator
is simply $A^\infty$, in both cases. $\ker\sigma$ does not contribute to the dimension
spectrum, due to prop.~\ref{prop:ker}.

It remains to consider an element of the form 
$\;f=\sum_{n\in\NN}(f_nw^n+f_{-n-1}(w^*)^{n+1})\,\in A^\infty$. 
All the terms in $f$ are off-diagonal, but the
one proportional to $f_0=\frac{1}{2\pi}\int_{S^1}\sigma(f)\de\theta$. Then
$\tr_{\ell^2(\N)}(fN^{-s})=f_0\,\zeta(s)$ and $\tr_{\HH_+}(f|D|^{-s})=2f_0\,\zeta(s-1)$,
where $\zeta(s)$ is the Riemann zeta-function. This concludes the proof.
\end{prova}

Notice that $1$ is not in the dimension spectrum of the second triple,
even though $S^1$ is a classical subspace of $S^2_0$, which seems quite curious.

%%% ======================================================================
\section{A spectral triple for $\mu(S^2_q)$.}\label{sec:mu}
Here we consider the spectral triple over $S^2_q$ associated to the 
representation $\mu = \mu_+\oplus\mu_-$ given by (\ref{eq:qdisk}). 
The Hilbert space is naturally isomorphic to $\ell^2(\N)\otimes\C^2$,
and the representation reads
\begin{equation*}
\mu(a)=w\sqrt{1-q^{4N}}\otimes\ma{1 & 0 \\ 0 & 1}\:,\qquad
\mu(b)=q^{2N}\otimes\ma{1 & 0 \\ 0 & -1}\;,
\end{equation*}
while as the Dirac operator $D'$ we take
\begin{equation*}
D'=N\otimes\ma{0 & 1 \\ 1 & 0}\:,
~ {\rm so}~
F:={\rm sign}\,D' = id_{\ell^2(\N)}\otimes\ma{0 & 1 \\ 1 & 0}
~ {\rm and}~
|D'|=N\otimes id_{\C^2}\;.
\end{equation*}
The resulting spectral triple can be thought of as a particular limit
of the one considered in~\cite[sec.~5]{CP1}%
\footnote{%
The generators $A,B$ used in~\cite{CP1} are related to ours generators
through $a=\lim_{c\to\infty}(c^{-1/2}B^*)$ and $b=q^2\lim_{c\to\infty}(c^{-1/2}A)$,
and the Dirac operator is shifted by $1$, in order to remove $0$ from the
spectrum. A finite shift of $D'$ results in a finite rank perturbation
of the sign $F$, and so it does not affect dimension spectrum, residues
and the index map.}.

Let $A^\infty$ be the algebra defined in sec.~\ref{sec:A}, and recall that
$\mu_{\pm}(a),\mu_{\pm}(b)\in A^\infty$. Then $\mu(a),\mu(b)\in A^\infty\oplus A^\infty$.
Since $\sigma\mu_{\pm}(a)=e^{i\theta},\sigma\mu_{\pm}(b)=0$, if we identify
$\A(S^2_q)$ with its representation $\mu$, then
\begin{equation}
\A(S^2_q)\subset\{(x,y)\in A^\infty\oplus A^\infty\,|\,\sigma(x)=\sigma(y)\}
=:C^\infty(S^2_q)\;,
\end{equation}
where $C^\infty(S^2_q)$ is the algebra anticipated in eq.~(\ref{eq:Cinf}).
From the direct sum of two copies of (\ref{eq:seq}), passing to the diagonal
($\sigma(x)=\sigma(y)$) we obtain the exact sequence:
\begin{equation}\label{eq:sseq}
0\to\rd\oplus\rd\to C^\infty(S^2_q)\stackrel{\rho}{\to} C^\infty(S^1)\to 0\;,
\end{equation}
where, as in the previous section, we call 
$\rho(x\oplus y):=\sigma(x)=\sigma(y)$, for $x\oplus y\in C^\infty(S^2_q)$.
Then, on the generators $a$ and $b$ the morphism $\rho$ is given by
$\rho(a)(\theta)=e^{i\theta}$, $\rho(b)(\theta)=0$.

Let $B$ be the algebra generated by $C^\infty(S^2_q)$ and commutators with $D'$.
Since $[D',\mu(a)]=\mu(a)F$ and $[D',\mu(b)]=-2N\mu(b)F$ are in $A^\infty
\otimes\mathrm{Mat}_2(\C)$ (since $Nq^{2N}\in A^\infty$), then:
\begin{equation*}
B\subset A^\infty\otimes\mathrm{Mat}_2(\C)\;.
\end{equation*}
Let $\delta':=[|D'|,\,.\,]$. From $\delta'\mu(a)=\mu(a)$ and $\delta'\mu(b)=0$
it follows that $B$ is $\delta'$-invariant, and hence it is in the smooth domain
of $\delta'$. We have thus proved the following proposition.
\begin{prop}
$(C^\infty(S^2_q),\ell^2(\N)\oplus\ell^2(\N),D')$ is a \emph{regular} spectral triple.
\end{prop}

\noindent 
Let, for $s\in\C$ with sufficiently large real part, 
the `zeta-type' function associated to $T\in B$ be given by
\begin{equation*}
\zeta_T(s):=\tr_{\ell^2(\N)\oplus\ell^2(\N)}(T|D'|^{-s})\;.
\end{equation*}
If we extend the $*$-algebra morphism $\sigma$ to:
\begin{equation}\label{eq:trho}
\tilde{\rho}=\sigma\otimes id:A^\infty\otimes\mathrm{Mat}_2(\C)
\to C^\infty(S^1)\otimes\mathrm{Mat}_2(\C)\:,
\end{equation}
then follows the Proposition \ref{2.4}.
\begin{prop}\label{2.4}
The dimension spectrum is $\{1\}$ and the residue is
\begin{equation*}
\mathrm{Res}_{s=1}\zeta_T(s)=\frac{1}{2\pi}\int_{S^1}\de\theta\;\tr_{\C^2}\tilde{\rho}(T)
\end{equation*}
for all $T\in B$, where $\tilde{\rho}$ is the map in (\ref{eq:trho}).
\end{prop}
\begin{prova}
Any $T\in B$ can be written as $T=x\otimes 1+y\otimes M$, with $x,y\in A^\infty$
and $M$ a traceless $2\times 2$ matrix.
Notice that $\zeta_{yM}(s)=0$ since $yM$ is off-diagonal due to the presence of $M$.

From the proof of the proposition \ref{prop:proof} we derive
\begin{equation*}
\zeta_T(s)=\zeta_x(s)=2\,\tr_{\ell^2(\N)}(xN^{-s})
   =\left(\frac{1}{2\pi}\int_{S^1}2\sigma(x)\de\theta\right)\zeta(s)\;.
\end{equation*}
To conclude the proof is sufficient to notice that from the definition
(\ref{eq:trho}) it follows: $\tr_{\C^2}\tilde{\rho}(T)\equiv 2\sigma(x)$.
\end{prova}

%%% ======================================================================
\subsection{Non-triviality of the Chern character}
Consider the projector 
$\ket{1}\!\left<1\right|\oplus 0$ over $\ell^2(\N)\oplus\ell^2(\N)$.
It has finite rank, hence it is an element of $\rd\oplus\rd$ and its image in
$C^\infty(S^2_q)$ via the morphism in (\ref{eq:sseq}) is a projector too.
This projector generates, together with the trivial one, the $K_0$-group
of $S^2_q$~\cite{Mas}. The index of the associated %(trivial) 
twisted Dirac operator (c.f.~App.~\ref{app:A}) is easily computed as
in~\cite[sec.~5]{CP1}, to be equal to $1$,
proving the non-triviality of the Chern character of the spectral triple.

Another interesting projector $p'\in\A(S^2_q)\otimes\mathrm{Mat}_2(\C)$ is given by
\begin{equation}\label{eq:pB}
p'=\frac{1}{2}\ma{1+b & a^* \\ a & 1-q^{-2}b}
\end{equation}
($p'$ is associated to the $SU_q(2)$-principal bundle with base space
$S^2_q$ in~\cite{BM}\footnote{Using the notations of~\cite{BM}, when $s=1$ 
(the condition corresponding to the equatorial Podle\'s sphere), 
the elements $(\xi,\eta,\zeta)$ of~\cite{BM}
correspond to $(a^*,-a,-b)$ and their projector is $e_{s=1}\equiv p'$.}).

\pagebreak

Note that $C(S^2_q)^{\oplus 2}p'$, completed with respect to 
a suitable inner product, is the analogue of the tautological line bundle.
In fact when $q=1$, $a\to 2z(1+z\bar{z})^{-1}$ and
$1-b\to 2(1+z\bar{z})^{-1}$, where $z=e^{i\varphi}\cot\theta/2$
is the stereographic coordinate on $S^2$. Thus when $q=1$, the projector $p'$
becomes the celebrated Bott projector
\begin{equation*}
p_B(z)=\frac{1}{1+z\bar{z}}\ma{z\bar{z} & \bar{z} \\ z & 1},\quad
p_B(\infty)=\ma{1 & 0 \\ 0 & 0}\;,
\end{equation*}
which, together with $[1]$ generates $K^0(S^2)$.\\
In the rest of this section we show that $\mathrm{ch}^F([p'])=-1$,
proving that $1-p'$ and $\ket{1}\!\left<1\right|\oplus 0$ are equivalent projectors
and that $p'$ can be taken as non-trivial generator of $K_0(C(S^2_q))$.

Denote $V_{\pm}\simeq\ell^2(\N)$ the two components of the Hilbert space
with the canonical orthonormal basis $\ket{n}_{\pm}$. 
To deal with $p'$ we need to lift the representation of the algebra 
and the sign of the Dirac operator to $2\times 2$ matrices, 
so let $V'_{\pm}=\C^2\otimes V_{\pm}$ and $F'=id_{\C^2}\otimes F$.

We choose the following (orthonormal) basis of $V'_\pm$
\begin{align*}
\ket{n}_\pm^0 &:=\sqrt{\frac{1-rq^{2n}}{2}}\bn{1}{0}\otimes\ket{n}_\pm
             -\sqrt{\frac{1+rq^{2n}}{2}}\bn{0}{1}\otimes\ket{n+1}_\pm\;,
             &\forall\;r=\pm,\;n\in\N\;, \\
\ket{n}_\pm^1 &:=\sqrt{\frac{1+rq^{2n}}{2}}\bn{1}{0}\otimes\ket{n}_\pm
             +\sqrt{\frac{1-rq^{2n}}{2}}\bn{0}{1}\otimes\ket{n+1}_\pm\;,
             &\forall\;r=\pm,\;n\in\N\;, \\
\ket{0}_{\pm}   &:=\bn{0}{1}\otimes\ket{1}_{\pm}\;.
\end{align*}
%From the equalities:
A straightforward calculation shows that
\begin{equation*}
\mu_\pm(p')\ket{n}^0_\pm=0\:,\qquad \mu_+(p')\ket{0}_+=0\:,\qquad
\mu_\pm(p')\ket{n}^1_\pm=\ket{n}^1_\pm\:,\qquad
\mu_-(p')\ket{0}_-=\ket{0}_-\:.
\end{equation*}
From these we deduce that $\{\ket{n}_+^1\}$ is a basis for $\mu_+(p')V'_+$
and $\{\ket{n}_-^1,\ket{0}_-\}$ is a basis for $\mu_-(p')V'_-$.

Let $F_{p'}:=\mu_-(p')F'\mu_+(p')$, it maps $\mu_+(p')V'_+$ to $\mu_-(p')V'_-$.
Then
\begin{align*}
F_{p'}\ket{n}^1_+ &=\sqrt{1-q^{4n}}\ket{n}^1_-\neq 0\:, &\forall\;n\in\Z_+\:, \\
F_{p'}^*\ket{n}^1_- &=\sqrt{1-q^{4n}}\ket{n}^1_+\neq 0\:, &\forall\;n\in\Z_+\:, \\
F_{p'}^*\ket{0}_- &=0\;.
\end{align*}
Hence  $\ker F_{p'}=0$, $\ker F_{p'}^*\simeq\C$ and 
$\mathrm{ch}^F([p'])=\mathrm{Index}(F_{p'})=-1$.

%%% ======================================================================
\subsection{Local index formula}
Consider the $*$-linear map $\mu_0:=\mu_+-\mu_-:C^\infty(S^2_q)\to\B(\ell^2(\N))$.
We claim that $\mu_0$ has image in $\rd$. 
In fact, the equality $\mu_0(xy)=\mu_0(x)\mu_+(y)+\mu_-(x)\mu_0(y)$ 
implies that $\mu_0(xy)\in\rd$ if $\mu_0(x),\mu_0(y)\in\rd$ 
($\rd$ being a two-sided $*$-ideal in $A^\infty$). 
Moreover the generators of the algebra satisfy $\mu_0(a)=0,\mu_0(b)=2q^{2N}\in\rd$.
As a consequence the commutator
\begin{equation*}
[F,\mu(x)]=\mu_0(x)\otimes\textrm{\footnotesize$\ma{0 & -1 \\ 1 & 0}$}
\end{equation*}
is traceclass for all $x\in C^\infty(S^2_q)$ and the map $\mathrm{ch}^F_0$ in
eq.~(\ref{eq:ch}) is well-defined. \pagebreak

The pairing between the cyclic cocycle defined by the map $\mathrm{ch}^F_0$
and the class $[p]\in K_0(\A)$ of a projector $p$ gives the index formula
\begin{equation*}
\mathrm{Index}(pFp)=\tfrac{1}{2}\tr_{\ell^2(\N)\oplus\ell^2(\N)}(\gamma F[F,\mu(p)])
=\tr_{\ell^2(\N)}\mu_0(p)\;,
\end{equation*}
where $\gamma = \textrm{\footnotesize$\ma{1 & 0 \\ 0 & -1}$}$ is the grading.
Theorem~\ref{thm}, applied to the present case, states that $\mathrm{ch}^F$ is
cohomologous to the cocycle with only one component $\varphi_0$ given by
\begin{equation*} 
\varphi_0(x)=\mathrm{Res}_{s=0}s^{-1}\tr_{\ell^2(\N)\oplus\ell^2(\N)}
\bigl(\gamma\mu(x)|D'|^{-2s}\bigr)
=\mathrm{Res}_{s=0}s^{-1}\tr_{\ell^2(\N)}\bigl(\mu_0(x)N^{-2s}\bigr)\;.
\end{equation*}
But $\psi(s):=\tr_{\ell^2(\N)}\bigl(\mu_0(x)N^{-2s}\bigr)$ is holomorphic since 
$\mu_0(x)\in\rd$. Thus, $\mathrm{Res}_{s=0}s^{-1}\psi(s)=\psi(0)$ 
and $\varphi_0(p)\equiv\mathrm{ch}^F_0([p])$.

Interestingly, theorem~\ref{thm} is here trivially satisfied 
(i.e.~the coboundary `$\mathrm{ch}^F-\varphi$' is zero). 
We shall encounter a similar situation in the next section.

%%% ======================================================================
\section{A spectral triple for $\pi(S^2_q)$.}\label{sec:pi}

In this section we discuss the spectral triple of~\cite{DLPS}.
To simplify the story, we work with polynomial algebra $\A(S^2_q)$.
Analogous results hold also for $C^\infty(S^2_q)$, but are more complicated.
The Hilbert space is $\HH=\HH_+\oplus\HH_-\simeq\HHH\otimes\C^2$, where
$\ket{l,m}$ is declared to be an orthonormal basis of $\HHH$, $m=-l,-l+1,...,l$.

We identify $\HH_{\pm}$ with $\HHH$ through the isometry $\ket{l,m}_{\pm}\to\ket{l,m}$
and consider the representations $\pi_{\pm}$ at page \pageref{eq:Drep} as representations
over $\HHH$. Then, if we call
\begin{equation*}
\rho_{\pm}:=(\pi_+\pm\pi_-)/2:\A(S^2_q)\to\B(\HHH)
\end{equation*}
the spinorial representation $\pi=\pi_+\oplus\pi_-$ can be rewritten as
\begin{equation*}
\pi=\rho_+\otimes id_{\C^2}+\rho_-\otimes\gamma\;,
\end{equation*}

\vspace{-5pt}

\noindent%
where $\gamma:=\,${\footnotesize$\ma{1 & 0 \\ 0 & -1}$} is the grading.
As stated in the introduction, the Dirac operator is $D=|D|\otimes F$, where
\begin{equation*}
|D|\ket{l,m}=(l+\tfrac{1}{2})\ket{l,m}\;,
\qquad\quad F=\textrm{\footnotesize$\ma{0 & 1 \\ 1 & 0}$}\;.
\end{equation*}
Let us rewrite explicitly $\rho_{\pm}(a)$ and $\rho_{\pm}(b)$:
\begin{subequations}\label{subeq}
\begin{align}
\rho_+(a)\ket{l,m} = 
&\;\frac{\sqrt{1-q^{2(l+m+1)}}\sqrt{1-q^{2(l+m+2)}}}{1-q^{4(l+1)}}\ket{l+1,m+1} \\
& -\frac{\sqrt{q^{2(l+m)}-q^{4l}}\sqrt{q^{2(l+m+1)}-q^{4l}}}{1-q^{4l}}\ket{l-1,m+1}\;, 
\label{subeq:b} \\
\rho_+(b)\ket{l,m} = 
&\,-\frac{\sqrt{1-q^{2(l+m+1)}}\sqrt{q^{2(l+m+2)}-q^{4l+6}}}{1-q^{4(l+1)}}\ket{l+1,m} \\
& -\frac{\sqrt{1-q^{2(l+m)}}\sqrt{q^{2(l+m+1)}-q^{4l+2}}}{1-q^{4l}}\ket{l-1,m}\;,\\
\rho_-(a)\ket{l,m} = 
&\,\frac{(1-q^4)q^{3l+m}}{(1-q^{2l})(1-q^{2l+2})}\sqrt{1-q^{2(l+m+1)}}
                       \sqrt{1-q^{2(l-m)}}\ket{l,m+1}\;, \label{subeq:e} \\
\rho_-(b)\ket{l,m} = 
&\,\frac{(1-q^2)q^{2l+1}}{(1-q^{2l})(1-q^{2l+2})}\Big\{1+q^{4l+2}-(1+q^2)
                       q^{2(l+m)}\Big\}\ket{l,m}\;. \label{subeq:f}
\end{align}
\end{subequations}
Notice that $\rho_-(a)$ and $\rho_-(b)$ do not act on the index $l$, while
$\rho_+(a)$ and $\rho_+(b)$ naturally decompose as the sum of two operators
that shift the index $l$ of $\pm 1$. With obvious notations, call $\rho_+(a)=a_++a_-$
and $\rho_+(b)=b_++b_-$. Recall that $\delta:=[|D|,\,.\,]$.
\begin{prop}
Let $B$ be the $*$-algebra generated by $\pi(\A(S^2_q))$
and $[D,\pi(x)]$, $x\in\A(S^2_q)$. Then $B\subset\B(\HH)$
and $(\A(S^2_q),\HH,D)$ is a $2^+$-dimensional regular spectral triple.
\end{prop}
\begin{prova}
$\delta(a_{\pm})=\pm a_{\pm}$ and $\delta(b_{\pm})=\pm b_{\pm}$
are bounded and in the smooth domain of $\delta$.
Since
\begin{equation*}
[D,\pi(a)]=(a_+-a_-)\otimes F+2\rho_-(a)D\otimes\gamma\qquad
[D,\pi(b)]=(b_+-b_-)\otimes F+2\rho_-(b)D\otimes\gamma\;\;,
\end{equation*}
then $B$ is in the $*$-algebra generated by $\{1,a_{\pm},b_{\pm},F\}$,
modulo terms in the kernel of $\delta$, and hence $B$ is in the smooth
domain of $\delta$. The metric dimension is $2$ by isospectrality.

The boundedness of the terms linear in $D$ is a consequence of the fact
that $\rho_-$ maps $\A(S^2_q)$ into rapid decay matrices. This will be
proved in the next proposition.
\end{prova}

\begin{prop}\label{prop1}
$\rho_-$ maps $\A(S^2_q)$ in $\rd$. Therefore, $[%id_{\HHH}\otimes 
F,\pi(x)]=2\rho_-(x)\otimes F\gamma$
is traceclass for all $x\in\A(S^2_q)$ and the associated Fredholm module is finite summable.
\end{prop}
\begin{prova}
Though $\rho_-$ is not a representation, it satisfies the following identity
\begin{equation*}
\rho_-(xy)=\pi_+(x)\rho_-(y)+\rho_-(x)\pi_-(y) \ .
\end{equation*}
Since $\rd$ is a two-sided $*$-ideal in $C^\infty(S^2_q)$, and $\A(S^2_q)\subset C^\infty(S^2_q)$,
then $\rho_-(xy)\in\rd$ if $\rho_-(x)$ and $\rho_-(y)\in\rd$.
So, to prove that $\rho_-$ has image in $\rd$ we need just to do
the check for the generators of the algebra.

In (\ref{subeq:e}), $q^{l+m}\leq 1$, the square roots are no greater than $1$
and $(1-q^\alpha)^{-1}\leq (1-q)^{-1}\;\forall\;\alpha\geq 1$. In (\ref{subeq:f})
the quantity in the big parenthesis is (in modulus) $\leq 1-q^{4l}\leq 1$.
Then, we deduce that the (nonzero) matrix coefficients of 
$\rho_-(a)$ and $\rho_-(b)$ satisfy
\begin{equation*}
|\inner{l,m+1|\rho_-(a)|l,m}| \leq (1-q)^{-2}q^{2l}\qquad\quad
|\inner{l,m|\rho_-(b)|l,m}| \leq (1-q)^{-2}q^{2l}
\end{equation*}
and so they are rapid decay matrices.
\end{prova}

\noindent
(The Fredholm module $(F, \pi)$ is 1-summable according to terminology of
\cite{Con0} and 0-summable according to \cite{Hig}.)

\nopagebreak[4]

%%% ======================================================================
\subsection{An approximate representation}
To simplify the computations, it is useful to cut smoothing contributions.

For all $\alpha>0$, $(1-q^{\alpha l})^{-1}-1$ is a rapid decay sequence.
For $0\leq u\leq 1$, $|1-\sqrt{1-u}|\leq u$.

\pagebreak

\noindent
Then,
\begin{equation*}
|q^{l+m}-\sqrt{q^{2(l+m)}-q^{4l}}|=q^{l+m}|1-\sqrt{1-q^{2(l-m)}}|\leq q^{3l-m}\leq q^{2l}
\end{equation*}
is a rapid decay sequence and the first square root in (\ref{subeq:b}) coincide with
$q^{l+m}$ modulo rapid decay sequences. Applying the same argument to the other
square roots in (\ref{subeq}), one prove that the operators $\lambda(a)$
and $\lambda(b)$, defined by
\begin{align*}
\lambda(a)\ket{l,m} &:=\sqrt{1-q^{2(l+m+1)}}\sqrt{1-q^{2(l+m+2)}}\ket{l+1,m+1}
      -q^{2(l+m)+1}\ket{l-1,m+1}\;\;, \\
\lambda(b)\ket{l,m} &:=-q^{l+m+2}\sqrt{1-q^{2(l+m+1)}}\ket{l+1,m}
      -q^{l+m+1}\sqrt{1-q^{2(l+m)}}\ket{l-1,m}\;\;,
\end{align*}
differ from $\pi_{\pm}(a)$ and $\pi_{\pm}(b)$ by a rapid decay matrix.

The closure in the operator norm of $\rd$ is the two-sided $*$-ideal $\K$.
If $\tilde{\lambda}$ is the projection into the Calkin algebra $\B(\HH)/\K$,
\begin{equation*}
0\to\K\to C(S^2_q)\stackrel{\tilde{\lambda}}{\to}C(S^2_q)/\K\to 0
\end{equation*}
then, $\lambda(a)$ and $\lambda(b)$ are representatives of $\tilde{\lambda}(a)$ and
$\tilde{\lambda}(b)$, and the $C^*$-algebra they generate coincide with $C(S^2_q)$
modulo compact operators.

When considering $\A(S^2_q)$, the algebra of polynomials in $\lambda(a)$ and $\lambda(b)$
coincide with $\A(S^2_q)$ modulo rapid decay matrices,
and the difference can be neglected when computing zeta-type functions.

The derivatives $\partial:=[D,\,.\,]$ and $\delta:=[|D|,\,.\,]$ send rapid decay matrices to rapid
decay matrices, and then the operators $[D,\pi(a)]-[D,\lambda(a)\otimes id_{\C^2}]$,
$[D,\pi(b)]-[D,\lambda(b)\otimes id_{\C^2}]$,
$\delta\bigl(\pi(a)\bigr)-\delta\bigl(\lambda(a)\bigr)$,
$\delta\bigl(\pi(b)\bigr)-\delta\bigl(\lambda(b)\bigr)$
are all rapid decay matrices and we have an ``approximate representation'' of
all the algebra $\bigcup_{k\in\NN}\delta^k(B)$ modulo $\rd$.

%%% ======================================================================
\subsection{The dimension spectrum}

Let $U\ket{l,m}=q^{l+m}\ket{l,m}$ and $\mathcal{Q}_q$ be the (two-sided $*$-)ideal
in $\B(\HH)$ generated by the two elements $U$ and $V=1-\sqrt{1-(qU)^2}$ (notice
that they are not compact, $||U||_E=1$ on each subspace $E\subset\HHH$ 
of finite codimension).\\
Any operator in $\mathcal{Q}_q$ can be written in the form 
$T=xUy$ or $T=xVy$ for some bounded $x, y$, and so
they satisfy (since $1-\sqrt{1-q^{2(l+m+1)}}\leq q^{2(l+m+1)}\leq q^{l+m}$)
\begin{equation*}
|\inner{l,m|T|l,m}|\leq ||x||\; ||y||\,q^{l+m}\;.
\end{equation*}
Then, for all $T\in\mathcal{Q}_q$,
\begin{align*}
|\zeta_T(s)| & \leq ||x||\; ||y||\sum\nolimits_{l+\frac{1}{2}\in\N}(l+\tfrac{1}{2})^{-s}
\sum\nolimits_{l+m=0}^{2l}q^{l+m} \\
&=||x||\; ||y||\tfrac{1}{1-q}\zeta(s)+\textrm{holomorphic function}\;.
\end{align*}
So $\zeta_T(s)$ has a unique residue in $s=1$. This is not identically zero, since
for $T=b^2$ we have
\begin{align*}
\zeta_{b^2}(s) &=\sum\nolimits_{l+\frac{1}{2}\in\N}(l+\tfrac{1}{2})^{-s}
    2q^2\sum\nolimits_{k=0}^{2l}(q^{2k}-q^{4k})
 \\ &=\tfrac{2q^{4}}{1-q^2}\zeta(s)+\textrm{holomorphic function}
\end{align*}
and then
\begin{equation*}
\mathrm{Res}_{s=1}\zeta_{b^2}(s)=\tfrac{2q^{4}}{1-q^2}\neq 0\;.
\end{equation*}
Let $(\rd\cup\mathcal{Q}_q)$ be the ideal generated by $\rd\cup\mathcal{Q}_q$.
Let $\nu:\A(S^2_q)\to \B (\HH)$ be the representation defined by
\begin{equation*}
\nu(a)\ket{l,m}=\ket{l+1,m+1}\;,\qquad \nu(b)=0\;.
\end{equation*}
Then for $x\in\A(S^2_q)$, $\pi(x)-\nu(x)\in(\rd\cup\mathcal{Q}_q)$,
or with a slight abuse of notation 
$\nu$ is a projection from $\A(S^2_q)$ to $\A(S^2_q)/(\rd\cup\mathcal{Q}_q)$.

The ideal $\mathcal{Q}_q$ is $\partial$-invariant and $\delta$-invariant.
Then, the computation of the dimension spectrum reduces to
the computation of the dimension spectrum of $\nu(\A(S^2_q))$,
i.e.~of $\A(S^2_0)$ in the isospectral representation (\ref{pi0}).

The contribution of $\mathcal{Q}_q$ gives a simple pole at $s=1$,
while the contribution of $\nu(\A(S^2_q))$ gives a simple pole at $s=2$
(Proposition \ref{2.4}).
Therefore we have
\begin{prop}
The dimension spectrum is $\Sigma=\{1,2\}$.
\end{prop}
We compute now the residue at $s=2$.
Let $x\in\A(S^2_q)$ and call $T=\nu(x)\otimes id_{\C^2}$. 
From Proposition~\ref{prop:proof} we have
\begin{equation*}
\mathrm{Res}_{s=2}\zeta_{T}(s)=
2\,\mathrm{Res}_{s=2}\tr_{\HHH}(T|D|^{-s})=\frac{2}{\pi}\int_{S^1}\sigma(T)\de\theta\;.
\end{equation*}
The last residue is equal to the residue of $\zeta_x(s)$.
Since $\sigma(\nu(a))=e^{i\theta}$ and $\sigma(\nu(b))=0$,
the final expression for the `noncommutative integral' of $x\in\A(S^2_q)$ is
\begin{equation}
\mathrm{Res}_{s=2}\zeta_x(s)=\frac{2}{\pi}\int_{S^1}\rho(x)\de\theta\;,
\end{equation}
where $\rho:\A(S^2_q)\to\A(S^1)$ is the $*$-algebra morphism defined
by $\rho(a)=e^{i\theta}$, $\rho(b)=0$.

%%% ======================================================================
\subsection{Non-triviality of the Chern character}
We show that the Chern character of the spectral triple is not trivial,
computing the pairing with the class of the projector $p'$ defined
in eq.~(\ref{eq:pB}).
From prop.~\ref{prop1}, the commutator
\begin{equation*}
[F,\pi(x)]=2\rho_-(x)\otimes\textrm{\footnotesize$\ma{0 & -1 \\ 1 & 0}$}
\end{equation*}
is traceclass for all $x\in\A(S^2_q)$ and the map 
$\mathrm{ch}^F_0$ in eq.~(\ref{eq:ch}) is well-defined,
\begin{equation*}
\mathrm{ch}^F_0(x)=2\tr_{\HHH}\rho_-(x)\;.
\end{equation*}
Then,
\begin{equation*}
\mathrm{ch}^F([p'])=(1-q^{-2})\,\tr_{\HHH}\rho_-(b)
    =-q^{-2}(1-q^2)^2\sum_{l,m}\frac{[l-m+1][l+m]}{[2l][2l+2]} \;\; .
\end{equation*}
Performing the sum in $m=-l,...,l$ and setting $n=l+\tfrac{1}{2}$, we obtain
\begin{equation*}
-\mathrm{ch}^F([p'])=\sum_{n=1}^\infty\frac{2n(1-x)^2(1+x^{2n})-(1-x^2)(1-x^{2n})}
{(1-x^{2n+1})(1-x^{2n-1})}x^{n-1}=:\sum_{n=1}^\infty f_n(x)=f(x)
\end{equation*}
where $x=q^2\in\,]0,1[$, and $f(x)$ is an integer-valued function (being the index
of a Fredholm operator) that we want to compute.

From the inequality $|f_n(x)|\leq(4n+2)x^{n-1}$ we deduce (Weierstrass M-test) that
the series is absolutely (hence uniformly) convergent in each interval $[0,x_0]
\subset [0,1[\,$. Then, it converges to a function $f(x)$ that is continuous in $[0,1[\,$.
A continuous function $f:\,\,]0,1[\,\to\Z$ is constant. By continuity, $f(x)$ is constant
in $[0,1[$ and can be computed setting $x=0$.
Since $f_n(0)=\delta_{n,1}$, we deduce that $f(x)=f(0)=1$ and so
\begin{equation*}
\mathrm{ch}^F([p'])=-1\;\;,\qquad\textrm{for all}\;0<q<1\,.
\end{equation*}

%%% ======================================================================
\subsection{Local index formula}

Thm.~\ref{thm}, applied to the present case, states that $\mathrm{ch}^F$ is
cohomologous to the cocycle with two components $(\varphi_0,\varphi_2)$, given by
\begin{align*} 
\varphi_0(a_0) & =\mathrm{Res}_{s=0}s^{-1}\tr(\gamma a_0|D|^{-2s}) \;\;,\\
\varphi_2(a_0,a_1,a_2) & =\mathrm{Res}_{s=0}\tr(\gamma a_0
[D,a_1][D,a_2]|D|^{-2(s+1)})\;\;.
\end{align*}
Neglecting rapid decay matrices, we have $[D,x]\sim\delta\rho_+(x)$ and then
$\varphi_2$ is identically zero
\begin{align*}
\varphi_2(a_0,a_1,a_2) &\equiv\mathrm{Res}_{s=0}\,\tr(
\rho_+(a_0)\,\delta\rho_+(a_1)\,\delta\rho_+(a_2)|D|^{-2(s+1)}\otimes\gamma)=0\;\;,
\end{align*}
since $\tr_{\C^2}\gamma=0$.

Moreover, $\psi(s)=\tr(\gamma a_0|D|^{-2s})=2\tr_{\HHH}\bigl(\rho_-(a_0)|D|^{-2s}\bigr)$
is holomorphic and
\begin{equation*}
\varphi_0(a_0):=\mathrm{Res}_{s=0}s^{-1}\psi(s)=\psi(0)\equiv\mathrm{ch}^F_0(a_0)\;.
\end{equation*}
As in Section~\ref{sec:mu}.2, also for the spinorial representation
and the isospectral Dirac operator the theorem~\ref{thm} is trivially satisfied 
(the coboundary `$\mathrm{ch}^F-(\varphi_0, \varphi_2)$' is zero), 
though the cocycle $\varphi$ could be expected to have higher dimensional components
(as the metric dimension is $2$).

%%% ======================================================================
\appendix
\section{Generalities about spectral triples}\label{app:A}
In this appendix we recall some material from~\cite{Con0,CM,Con1}
(see also \cite{NCG} and \cite{Hig} for a comprehensive exposition).\\
Let $(A,\HH,D)$ be a spectral triple and consider the following (unbounded)
derivations on $\B(\HH)$:
\begin{equation*}
\partial,\delta:\B(\HH)\to\B(\HH)\;,\qquad\partial\,T=[D,T]\;,\qquad\delta(T)=[|D|,T]
\end{equation*}
The triple is \emph{regular} if
\begin{equation*}
A\cup \partial A\subset\bigcap_{j\in\NN}\mathrm{dom}\,\delta^j\;.
\end{equation*}
The  operator $D$ is \emph{finite dimensional} if there exist a $d\in\R^+$
(called the \emph{metric dimension}) such that the singular values of $|D|^{-1}$
(assume that $D$ is invertible) are of order $n^{-1/d}$ when $n\to\infty$.\\
~\\
The `zeta-type' function 
\begin{equation*}
\zeta_T(s):=\tr_{\HH}\bigl(T|D|^{-s}\bigr)
\end{equation*}
associated to $T\in V=\bigcup_{j\in\NN}\delta^j(A\cup\partial A)$ is defined 
is defined (and holomorphic) for $s\in\C$ with $\mathrm{Re}\,s>d$.
% (Weierstrass theorem).

For a finite-dimensional spectral triple, it makes sense the following definition:
\begin{df}
A spectral triple has \emph{dimension spectrum} $\Sigma$ if{}f $\Sigma\subset\C$
is a countable set and for all $T\in V$, $\zeta_T(s)$ extends to a meromorphic function
on $\C$ with poles in $\Sigma$ as unique singularities.
\end{df}

Residues of zeta-type functions are traces on the algebra, and are used to
compute the pairing between $K$-theory and cyclic-cohomology.\\
Denoting $F=\mathrm{sign}\,D$, the bounded commutator condition implies
that $[F,x]$ is compact for all $x\in A$, and then $(A,\HH,F)$ is a Fredholm
module. Recall that a Fredholm module is \emph{finite summable} if for $k$
sufficiently large, $[F,a_0]\ldots[F,a_k]$ is traceclass for all $a_j\in A$.

From now on, we suppose the triple is even and denote by $\gamma$ the grading.
If the associated Fredholm module is finite summable, the map
\begin{equation}\label{eq:ch}
\mathrm{ch}^F_n(a_0,\ldots,a_n) =\frac{\Gamma(\tfrac{n}{2}+1)}{2n!}\tr(\gamma F[F,a_0]\ldots[F,a_n])
\end{equation}
defines a periodic cyclic cocycle $\mathrm{ch}^F$ (with all components equal to zero, 
but the $n$th that is $\mathrm{ch}^F_n$) whose periodic cyclic cohomology class 
is independent of $n$, 
for all $n$ even and sufficiently large. The pairing of $\mathrm{ch}^F$ with
$K$-theory gives the index map.

If $\phi=(\phi_0,\phi_2,...)$ is an element of the periodic cyclic cohomology
group PHC$^{ev}$ (only a finite number of components are different from zero)
and $p\in\mathrm{Mat}_\infty(A)$ a projector, the explicit formula for
the pairing is
\begin{equation*}
\inner{\phi,[p]}:=
\phi_0(p)+\sum_{k\in\N}(-1)^k\frac{(2k)!}{k!}\phi_{2k}(p-\tfrac{1}{2},p,\ldots,p)
\end{equation*}
and the index formula states that $\inner{\smash[b]{\mathrm{ch}^F}\!,[p]}$ is the index
of the twisted Dirac operator $pDp$ (or equivalently, $pFp$). In the case $\phi=\mathrm{ch}^F$,
there is only one non-zero component.

There exists a general theorem relating the index to residues of zeta-type
functions, which we quote in the case we are interested in:
\begin{thm}[Connes-Moscovici~\protect{\cite[thm.~II.3]{CM}}]\label{thm}
Let $(A,\HH,D)$ be a regular, even spectral triple (with finite metric dimension
$d$), with dimension spectrum $\Sigma$ made of simple poles.
Then, the following formulas define a $(b,B)$-cocycle with the same cyclic cohomology
class of the Chern character $\mathrm{ch}^F$ ($n$ even $\leq d$):
\begin{align*}
\varphi_0(a_0) &=\mathrm{Res}_{s=0}s^{-1}\tr(\gamma a_0|D|^{-2s}) \;\;,\\
\varphi_n(a_0,...,a_n) &=\sum_{k\in\NN^{\times n}}\frac{(-1)^k}{k_1!...k_n!}\alpha_k
\mathrm{Res}_{s=0}\tr\Big(\gamma a_0[D,a_1]^{(k_1)}...[D,a_n]^{(k_n)}|D|^{-2(|k|+\tfrac{n}{2}+s)}\Big)\;\;.
\end{align*}
Here $T^{(j+1)}=[D^2,T^{(j)}]\;\forall\;j\in\NN$, $T^{(0)}=T$
and $\alpha_k^{-1}=(k_1+1)(k_1+k_2+2)...(k_1+...+k_n+n)$.
\end{thm}
If $T$ is an order-zero operator, $T^{(j)}$ is of order $j$. Then, each residue
of the previous theorem can be written as a residue in $s=|k|+n$ of the zeta-type
function associated to a suitable order-zero operator. By definition, it is
zero if $|k|+n\notin\Sigma$. The finiteness of the metric dimension
guarantee that all terms in the sum with $|k|+n>d$ are zero, hence
it is a finite sum. Moreover, if $0\notin\Sigma$, then
$\varphi_0(a_0)=\tr(\gamma a_0|D|^{-2s})\big|_{s=0}$.\\
All the terms $\varphi_n$ with $n>0$ are local (i.e.~they don't care of traceclass
contributions). The only non-local term is $\varphi_0$.

The theorem is particularly interesting when $A$ is a pre-$C^*$-algebra, 
since in this case $A$ has the same $K$-theory as its $C^*$-algebra completion.

%%% ======================================================================
\section{The projective module $L^2(S^2_q)p$}\label{app:B}
Let us return to the projector $p$ in (\ref{eq:p}). The range
of $p$ is in the kernel of $a^*$: $a^*p=a^*-a^*=0$.
Notice also that $b^2p=q^4(1-aa^*)p=q^4p$.

Let us construct the $C(S^2_q)$-module $C(S^2_q)p$. 
A linear basis is $\{a^nb^mp, (a^*)^{n+1}b^np,\;n,m\in\NN\}$, 
but $(a^*)^{n+1}b^np\propto b^{n}(a^*)^{n+1}p=0$.
Moreover $a^nb^mp$ is proportional to $a^np$ or $a^nbp$ depending on the parity
of $m$ (since $b^{2n}p=q^{4n}p$). Then a minimal linear basis can be taken as
\begin{equation*}
\ket{n}_0=a^{n-1}p\;, \qquad\quad\ket{n_1}=a^{n-1}bp\;,
\end{equation*}
with $n\in\N$. It is immediate to compute:
\begin{equation*}
\left\{\begin{array}{l}
a\ket{n}_0=\ket{n+1}_0 \\
a\ket{n}_1=\ket{n+1}_1
\end{array}\right.\;,\qquad
\left\{\begin{array}{l}
a^*\ket{n}_0=(1-q^{4(n-1)})\ket{n-1}_0 \\
a^*\ket{n}_1=(1-q^{4(n-1)})\ket{n-1}_1
\end{array}\right.\;,\qquad
\left\{\begin{array}{l}
b\ket{n}_0=q^{2(n-1)}\ket{n}_1 \\
b\ket{n}_1=q^{2(n+1)}\ket{n}_0
\end{array}\right.\;.
\end{equation*}
We compute the inner product imposing that $a$ and $a^*$ are hermitian conjugates. 
It easily follows that $\ket{n}_s$ are orthogonal. We compute their
norm imposing (omitting the subscript $0,1$):
\begin{equation*}
\delta_{m,n+1}\inner{n+1|n+1}=\inner{m|a|n}=(a^*\ket{m})^\dag\ket{n}=
(1-q^{4n})\delta_{m,n+1}\inner{n|n}\;.
\end{equation*}
So $c_n=\inner{n|n}$ satisfies $c_{n+1}=(1-q^{4n})c_n$ and the (positive)
solution is
\begin{equation*}
c_{n+1}^{0,1}=(1-q^{4n})(1-q^{4(n-1)})\ldots(1-q^4)c_0^{0,1}\;.
\end{equation*}
Now we can complete the space in this norm and obtain an Hilbert space.

We define the basis:
\begin{equation*}
\ket{n}_{\pm}=\frac{1}{\sqrt{2c_n^0}}\ket{n}_0\pm\frac{1}{\sqrt{2c_n^1}}\ket{n}_1\;, 
\qquad c^{0,1}_0\in\R^+\;.
\end{equation*}
If we fix $c_0^0=1$ and $c^1_0=q^4$, in this (orthonormal) basis
\begin{equation*}
a\ket{n}_\pm=\sqrt{1-q^{4n}}\ket{n+1}_{\pm}\;,\qquad
a^*\ket{n}_\pm=\sqrt{1-q^{4(n-1)}}\ket{n-1}_{\pm}\;,\qquad
b\ket{n}_\pm=\pm q^{2n}\ket{n}_{\pm}\;.
\end{equation*}
Then, the module associated with $p$ is equivalent to the representation
$\mu=\mu_+\oplus\mu_-$ over $\ell^2(\Z)\oplus\ell^2(\Z)$, discussed in
section \ref{sec1}.
It has not a $q\to 1$ limit (there are no scalar projectors in $C(S^2)$).

% ---- Bibliografia ------------------------------------------------------

%%% ======================================================================
\end{document}